\documentclass[12pt]{amsart}
\usepackage[margin=1in]{geometry}
\usepackage{microtype}
\usepackage[htt]{hyphenat}
\usepackage[english]{babel}
\usepackage[hyperindex, linktocpage]{hyperref}
\usepackage[alphabetic, backrefs, msc-links]{amsrefs}
\usepackage{amsmath,amssymb,amsfonts,amsthm,enumerate}
\usepackage{url}
\usepackage{graphicx}
\usepackage{xcolor}
\usepackage{esint}

\usepackage[T1]{fontenc}

\usepackage[cal=boondox]{mathalfa}

\definecolor{darkblue}{RGB}{0,0,170}
\definecolor{brickred}{RGB}{200,0,0}
\hypersetup{colorlinks=true, linkcolor=darkblue, citecolor=blue} 

\newtheorem{theorem}{Theorem}

\newcommand{\R}{\mathbb{R}}

\newcommand{\N}{\mathbb{N}}

\usepackage{yfonts}

\renewcommand{\varpi}{\omega}

\renewcommand{\le}{\leqslant}

\renewcommand{\ge}{\geqslant}

\renewcommand{\epsilon}{\varepsilon}

\newcommand{\e}{\varepsilon}

\newcommand{\PV}{\,{\rm P.V.}\,}

\numberwithin{equation}{section} 

\title[All functions are $s$-harmonic]{All functions are (locally) $s$-harmonic\\
(up to a small error) -- and applications}\thanks{Supported by th
Istituto Nazionale di Alta Matematica and the Australian Research Council Discovery Project
{\em N.E.W.} ``Nonlocal Equations at Work''.}

\author{Enrico Valdinoci}

\address{{\em Enrico Valdinoci:} 
Dipartimento di Matematica, Universit\`a degli studi di Milano,
Via Saldini 50, 20133 Milan, Italy, and
Istituto di Matematica Applicata e Tecnologie Informatiche,
Via Ferrata 1, 27100 Pavia, Italy,
and School of Mathematics and Statistics,
University of Melbourne, 813 
Swanston St, Parkville VIC 3010, Australia.}

\email{enrico@mat.uniroma3.it}


\begin{document}

\begin{abstract}
The classical and the fractional Laplacians exhibit a number of similarities,
but also some rather striking, and sometimes surprising, structural differences.

A quite important example of these differences
is that any function (regardless of its shape) can be locally approximated
by functions with locally vanishing fractional Laplacian,
as it was recently proved by Serena Dipierro, Ovidiu Savin and myself.

This informal note is an exposition of this result
and of some of its consequences.
\end{abstract}

\maketitle

\section{Introduction}

Given~$s\in(0,1)$, we take into account the so-called $s$-fractional Laplacian
\begin{equation}\label{D}
(-\Delta)^s u(x):=\int_{\R^n} \frac{2u(x)-u(x+y)-u(x-y)}{|y|^{n+2s}}\,dy.\end{equation}
In this definition, $u$ is supposed to be a sufficiently smooth function (to make the
integral convergent for small~$y$) and with some growth control at infinity
(to make the integral convergent for large~$y$).
Also, for the sake of simplicity,
a normalizing constant is dropped in~\eqref{D}. It is also interesting to
observe that, by splitting two integrals
and changing variables, equation~\eqref{D}
can be written as
\begin{equation}\label{D2}\begin{split}
(-\Delta)^s u(x)\,&=\lim_{\rho\searrow0}
\int_{\R^n\setminus B_\rho} \frac{u(x)-u(x+y)}{|y|^{n+2s}}\,dy
+\int_{\R^n\setminus B_\rho} \frac{u(x)-u(x-y)}{|y|^{n+2s}}\,dy
\\
&=2\lim_{\rho\searrow0}
\int_{\R^n\setminus B_\rho} \frac{u(x)-u(x+y)}{|y|^{n+2s}}\,dy\\
&=2\lim_{\rho\searrow0}
\int_{\R^n\setminus B_\rho(x)} \frac{u(x)-u(y)}{|x-y|^{n+2s}}\,dy\\
&=:2\PV
\int_{\R^n} \frac{u(x)-u(y)}{|x-y|^{n+2s}}\,dy,\end{split}
\end{equation}
where the notation ``$\PV$'' stands for ``in the Cauchy Principal Value Sense
(and the factor~$2$ will not be relevant for our purposes).

The fractional Laplacian is one of the most widely studied operators in the recent literature,
probably in view of its intrinsic beauty (in spite of the first impression that the definition in~\eqref{D}
can produce), of the large variety of different problems related to it,
and of its great potentials in modeling real-world phenomena
in applied sciences.

The setting in~\eqref{D} is clearly related to an ``incremental quotient'' of~$u$
which gets averaged in all the space. Indeed, roughly speaking, equation~\eqref{D} combines
together several special features related to the classical Laplacian:
\begin{enumerate}
\item[{\bf(I)}] The classical Laplacian arises
from a second order incremental quotient,
namely, for a smooth function~$u$ and a small increment~$h$, denoting by~$\{e_j\}_{j=1,\dots,n}$
the standard Euclidean basis of~$\R^n$, it holds that
\begin{eqnarray*}&&
2u(x)-u(x+he_j)-u(x-he_j)\\&=&2u(x)-\left(
u(x)+\nabla u(x)\cdot(he_j)+\frac12\,D^2u(x) (he_j)\cdot(he_j)+o(h^2)\right)\\&&\quad
-\left(
u(x)+\nabla u(x)\cdot(-he_j)+\frac12\,D^2u(x) (he_j)\cdot(he_j)+o(h^2)\right)
\\ &=& -h^2\partial^2_{jj} u(x)+o(h^2)\end{eqnarray*}
and so
$$ \lim_{h\to0} \frac{ 2u(x)-u(x+he_j)-u(x-he_j) }{h^2}=-\Delta u(x).$$
Comparing this with~\eqref{D}, we recognize
a structure related to incremental quotients in the definition of fractional Laplacian;
\item[{\bf(II)}] The classical Laplacian compares the value of a function with
its average. Indeed, for a small~$\rho>0$,
\begin{equation}\label{E0}\begin{split}
\int_{B_\rho(x)} u(y)\,dy \,&= \int_{B_\rho} u(x+y)\,dy \\&=
\int_{B_\rho}\left( u(x)+\nabla u(x)\cdot y+\frac12\,D^2u(x) y\cdot y+o(|y|^2)\right)\,dy
.\end{split}\end{equation}
Also, by odd symmetry we see that
$$ \int_{B_\rho} y_j\,dy=0\qquad{\mbox{ for all }}j\in\{1,\dots,n\}$$
and
$$ \int_{B_\rho} y_j\,y_k\,dy=0\qquad{\mbox{ for all }}j\ne k\in\{1,\dots,n\}.$$
Consequently, we can write~\eqref{E0} as
\begin{equation*}\begin{split}
\fint_{B_\rho(x)} u(y)\,dy \,&= 
u(x)+\frac12\,\sum_{j=1}^n\partial^2_{jj}u(x)\,
\fint_{B_\rho} y_j^2\,dy+o(\rho^2)\\
&= 
u(x)+\frac1{2n}\,\sum_{j=1}^n\partial^2_{jj}u(x)\,
\fint_{B_\rho} |y|^2\,dy+o(\rho^2)\\
&= 
u(x)+\frac{\rho^2}{2(n+2)}\,\Delta u(x)+o(\rho^2)
\end{split}\end{equation*}
and therefore
\begin{equation}\label{IM} 
-\Delta u(x)=2(n+2)\,\lim_{\rho\searrow0}
\fint_{B_\rho(x)} \frac{u(x)-u(y)}{\rho^2}\,dy.\end{equation}
Similarly,
\begin{equation*}\begin{split}
\fint_{\partial B_\rho(x)} u(y)\,d{\mathcal{H}}^{n-1}(y) \,&= 
\fint_{\partial B_\rho} 
\left( u(x)+\nabla u(x)\cdot y+\frac12\,D^2u(x) y\cdot y+o(|y|^2)\right)
\,d{\mathcal{H}}^{n-1}(y)
\\
&= 
u(x)+\frac12\,\sum_{j=1}^n\partial^2_{jj}u(x)\,\fint_{\partial B_\rho} 
y_j^2\,d{\mathcal{H}}^{n-1}(y)+o(\rho^2)
\\ &= u(x)+\frac1{2n}\,\sum_{j=1}^n\partial^2_{jj}u(x)\,\fint_{\partial B_\rho} 
|y|^2\,d{\mathcal{H}}^{n-1}(y)+o(\rho^2)\\&=
u(x)+\frac{\rho^2}{2n}\,\Delta u(x)+o(\rho^2)
\end{split}\end{equation*}
and therefore
\begin{equation}\label{IM2} \begin{split}-\Delta u(x)\,&
=2n\,\lim_{\rho\searrow0}
\fint_{\partial B_\rho(x)} \frac{u(x)-u(y)}{\rho^2}\,d{\mathcal{H}}^{n-1}(y)\\&
=2n\,\lim_{\rho\searrow0}
\fint_{\partial B_\rho(x)} \frac{u(x)-u(y)}{|x-y|^2}\,d{\mathcal{H}}^{n-1}(y)
.\end{split}\end{equation}
Once again, the factors~$2(n+2)$ and~$2n$
in~\eqref{IM} and~\eqref{IM2} are not important for our purposes,
but the similarities between~\eqref{D2}, \eqref{IM} and~\eqref{IM2} are evident
and suggest that the fractional Laplacian
is a suitably weighted average distributed
in the whole of the space.
\item[{\bf(III)}] The classical Laplace operator is variational and stems from a Dirichlet energy of the form
\begin{equation}\label{DOM1} \int |\nabla u(x)|^2\,dx.\end{equation}
Similarly, the fractional Laplacian is variational and the corresponding energy is the
Gagliardo-Slobodeckij-Sobolev seminorm
\begin{equation}\label{DOM2} \iint \frac{|u(x)-u(y)|^2}{|x-y|^{n+2s}}\,dx\,dy.\end{equation}
The integral in~\eqref{DOM1} usually ranges in a ``domain''~$\Omega\subseteq\R^n$
which should be considered as the region of space where ``action
takes place'', or, better to say the complement of the region
in which no action takes place (that is, the domain~$\Omega$
is the complement of the region~$\R^n\setminus\Omega$, where
the data of~$u$ are fixed). The fractional
counterpart of this is to take as ``natural domain'' for~\eqref{DOM2}
the complement (in~$\R^{2n}$) of the set~$(\R^n\setminus\Omega)
\times(\R^n\setminus\Omega)$ where the data of~$u(x)-u(y)$
are fixed, that is, it is common to integrate~\eqref{DOM2}
over the ``cross domain''
\begin{eqnarray*} Q_\Omega&:=&(\Omega\times\Omega)\cup
\big(\Omega\times(\R^n\setminus\Omega)\big)
\cup
\big((\R^n\setminus\Omega)\times\Omega\big)
\\&=&\R^{2n}\setminus
\big(
(\R^n\setminus\Omega)\times(\R^n\setminus\Omega)\big).
\end{eqnarray*}
\item[{\bf(IV)}] Most importantly, the fractional Laplacian enjoys ``elliptic'' features
that are similar to the ones of the classical Laplacian, e.g. in terms of maximum principle.
The regularizing effects of the fractional Laplacian can be somewhat ``guessed''
from the singularity of the integral kernel in~\eqref{D}: indeed, on the one hand, to
make sense of the integral in~\eqref{D}, one needs the function~$u$ to be ``smooth enough''
near~$x$; on the other hand, and somehow conversely, if the integral in~\eqref{D} is finite,
the function~$u$ needs to have some regularity property near~$x$, in order to compensate
the singularity of the kernel.
\end{enumerate}

Several classical and recent publications presented
the fractional Laplacian from different perspectives.
See in particular~\cites{MR0290095, MR0350027, MR2707618, MR2944369, MR3469920}.
In our postmodern world some
excellent online expositions of this topic
have also become available, see in particular the very useful webpage\\
{\tt https://www.ma.utexas.edu/mediawiki/index.php/Fractional\_Laplacian}

\medskip

We also recall that the fractional Laplacian can also be framed into the context of
probability and harmonic analysis, thus leading to different possible approaches
and several possible definitions, see~\cite{MR3613319}, and it is also
possible to provide a suitable setting in order to define the fractional Laplacian
for functions with polynomial growth at infinity, see~\cite{POLYN}.
\medskip

In spite of the extremely important similarities between the classical and the fractional
Laplacian, several structural differences between these operators arise.
See e.g.~\cite{BARI} for a collection of some of these basic differences.
Some of these differences have also extremely deep consequences
on some recent results in the theory of nonlocal equations, see\\
{\tt https://www.ma.utexas.edu/mediawiki/index.php/List\_of\_results\_that\_are\_fundamentally\_different\_to\_the\-local\_case}
\medskip

In this note, we recall one of the basic differences between the classical and the fractional
Laplacian, which has been recently discovered in~\cite{MR3626547}
and which presents a source of interesting consequences.
This difference deals with the so called ``$s$-harmonic functions'',
which are the (rather surprising) counterpart
of classical harmonic functions. 

The parallelism between
classical harmonic functions and $s$-harmonic functions lies in their definition,
since~$u$ is said to be harmonic (respectively, $s$-harmonic) at~$x$
if~$-\Delta u(x)=0$ (respectively, if~$(-\Delta)^s u(x)=0$).

Already from the definition, a basic difference between the classical
and the fractional case arises, since the definition of harmonic function at~$x$
only requires the function to be defined in an arbitrarily small neighborhood
of~$x$, while the 
definition of $s$-harmonic function requires the function to be globally
defined in~$\R^n$. This difference, which is somehow the counterpart
of the structural differences between~\eqref{D2} on one side
and~\eqref{IM} and~\eqref{IM2} on the other side, turns out to be
perhaps deeper than what may look
at a first glance. As a matter of fact,
the classical Laplacian is a very ``rigid'' operator, and for a function
to be harmonic some very restrictive geometric conditions must hold
(in particular, harmonic functions cannot have local minima).
In sharp contrast with this fact, the fractional Laplacian is very flexible
and the ``oscillations of a function that come from far''
can locally produce very significant contributions.

Probably, the most striking example of this phenomenon is that such
far-away oscillations can make the fractional Laplacian of any function
to almost vanish at a point, and in fact any given function,
without any restriction on its geometric properties, can be approximated
arbitrarily well by an $s$-harmonic function. In this sense, we have:

\begin{theorem}[``All functions are locally $s$-harmonic up to a small error''
\cite{MR3626547}]\label{TH}
For any~$\e>0$ and any 
function~$\bar v\in C^2(\overline{B_1})$, there exists~$v_\e$ such that
\begin{equation*}
\left\{\begin{matrix}
\|\bar v-v_\e\|_{C^2(B_1)}\le\e,\\
(-\Delta)^s v_\e=0{ \mbox{ in }}B_1.\end{matrix}\right.\end{equation*}\end{theorem}

A proof of this fact (in dimension~$1$ for the
sake of simplicity) will be given in Section~\ref{ALL:FUN}
(see the original paper~\cite{MR3626547} for the full
details of the argument in any dimension).\medskip

We stress that the phenomenon
described in Theorem~\ref{TH}
is very general, and it arises also for
other nonlocal operators, independently
from their possibly
``elliptic'' structure
(for instance all functions are locally $s$-caloric,
or $s$-hyperbolic, etc.), see~\cite{SCALOR}.

It is interesting to remark that the proofs
in~\cites{MR3626547, SCALOR} are not ``quantitative'',
in the sense that they are based on a contradiction argument,
and the ``shape'' of the approximating $s$-harmonic (or $s$-caloric, or $s$-hyperbolic)
function
cannot be detected by our methods. On the other hand, for a very nice quantitative
version of Theorem~\ref{TH} see 
Theorem~1.4 in~\cite{2017arXiv170806294R}.
See also~\cite{SALO}
for a quantitative approach to the parabolic case and~\cite{2016arXiv160909248G}
for related results (and, of course, quantitative proofs are harder
and technically more advanced than the one
that we present here).
In addition, 
results similar to Theorem~\ref{TH} hold true for nonlocal
operators with memory, see~\cite{ESAIM}.\medskip

Theorem~\ref{TH} possesses some simple, but quite interesting
consequences. In the forthcoming Section~\ref{APPE}
we present a few of them, related
to
\begin{enumerate}
\item[{\bf(i)}] The fractional Maximum Principle and Harnack Inequality;
\item[{\bf(ii)}] The classification of stable solutions for fractional equations;
\item[{\bf(iii)}] The diffusive strategy of biological populations.
\end{enumerate}

\section{Applications of Theorem~\ref{TH}}\label{APPE}

\subsection{The fractional Maximum Principle and Harnack Inequality}

One of the main features of the classical Laplace operator is that
it enjoys the Maximum Principle. For instance, as well known, it holds that:

\begin{theorem}\label{IA}
Let $u$ be a harmonic and nonnegative function 
in~$B_1$. If~$u(x_0)=0$ for some~$x_0\in B_1$,
then~$u$ is necessarily constantly equal to~$0$ in~$B_1$.
\end{theorem}

A classical quantitative version of Theorem~\ref{IA} was given by
Axel von Harnack
and can be stated as follows:

\begin{theorem}\label{HA1}
If~$u$ is 
harmonic in~$ B_1$ and nonnegative
in~$ B_1$, then, for every~$r\in(0,1)$,
$$ \sup_{B_r}u\le C_r\inf_{B_r} u ,$$
for some~$C_r>0$ depending on~$n$ and~$r$.
\end{theorem}

The original
manuscript by von Harnack is
available at \\
\begin{sloppypar}
{\tt https://ia902306.us.archive.org/9/items/vorlesunganwend\allowbreak00weierich/\allowbreak vorlesunganwend\allowbreak00weierich.pdf}
\end{sloppypar}
\medskip

The fractional counterpart of Theorem~\ref{HA1} goes as follows:

\begin{theorem}\label{HA2}
If~$u$ is 
$s$-harmonic in~$ B_1$ and nonnegative
in the whole of~$ \R^n$, then, for every~$r\in(0,1)$,
$$ \sup_{B_r}u\le C_r\inf_{B_r} u ,$$
for some~$C_r>0$ depending on~$n$ and~$r$.
\end{theorem}

See~\cites{MR1941020, MR1918242, MR2494809, MR2754080, MR2817382, MR3283558}
and the references therein for a detailed study of the fractional Harnack
Inequality. Of course, an important structural difference
between Theorems~\ref{HA1} and~\ref{HA2} (besides
the $s$-harmonicity versus the classical harmonicity) is the fact
that in Theorem~\ref{HA2} one requires a global condition on the
sign of the solution. Interestingly,
if in Theorem~\ref{HA2} one replaces the assumption that~$u$
is nonnegative in the full space with the assumption
that~$u$ is nonnegative just in the unit ball,
then the result turns out to be false, as described by the following
example:

\begin{theorem} \label{UAJ}
There exists a bounded function~$u$ 
which is $s$-harmonic in~$ B_1$, nonnegative
in~$ B_1$, not identically~$0$ in~$B_1$, but such that 
$$ \inf_{B_{1/2}} u= 0.$$\end{theorem}

Theorem~\ref{UAJ} suggests that some care has to be taken
when dealing with Maximum Principles and oscillation results
in the fractional case, and in fact the nonlocal character of the operator
requires global conditions for this type of results to hold,
in virtue of the contributions ``coming from far away''.\medskip

A proof of Theorem~\ref{UAJ} can be obtained directly from
Theorem~\ref{TH}. Indeed, we take~$n=1$, $\bar v(x):=x^2$
and~$\e:=\frac1{16}$. Then, Theorem~\ref{TH} provides a function~$v$
which is $s$-harmonic in~$(-1,1)$ and such that
$$ \|\bar v-v\|_{L^\infty((-1,1))}\le
\|\bar v-v\|_{C^2((-1,1))}\le\frac1{16}.$$
In particular, if~$|x|\ge\frac12$,
$$ v(x)\ge \bar v(x)-\frac1{16} =|x|^2-\frac1{16}
\ge\left( \frac12\right)^2-\frac1{16}=
\frac{3}{16},$$
while
$$ v(0)\le \bar v(0)+\frac1{16}=\frac1{16}.$$
Accordingly,
$$ \inf_{(-1,1)} v\le \frac1{16}<\frac3{16}\le
\inf_{(-1,1)\setminus(-1/2,1/2)} v,$$
which gives that
$$ \inf_{(-1,1)} v=\min_{[-1/2,1/2]} v=:\iota.$$
Then the function~$u:=v-\iota$ satisfies the thesis of
Theorem~\ref{UAJ}, as desired.
\medskip

For different approaches to the counterexamples to the local
Harnack
Inequality in the fractional setting see~\cite{MR2137058, FAILS} and also\footnote{We take
this opportunity to amend a typo in Theorem~2.3.1
of~\cite{MR3469920}, where~$\inf_{B_{1}} u$ has to
be replaced by~$\inf_{B_{1/2}} u$.}
Chapter~2.3 in~\cite{MR3469920}.

\subsection{The classification of stable solutions for fractional equations}

In the Calculus of Variations\footnote{Notice that the notion
of ``stability'' differs
from one scientific community
to another. In particular, the notion of stability
that we treat here does not agree with that in Dynamical Systems
or  Algebraic Geometry.}
literature, a solution~$u$ is called ``stable''
if it is the critical point of an energy functional whose
second variation is nonnegative definite at~$u$.
For instance, local minimizers of the energy are stable solutions,
and it is in fact often convenient to study stable
solutions since the stability class
is often preserved under suitable limit procedures and
it is sometimes technically easier
(or at least less difficult) to prove that a solution is stable
rather than deciding whether or not it is minimal.

We refer to the very nice monograph~\cite{MR2779463}
for a throughout discussion of the notion of stability and for many
related results.

A classical result in the framework of stable solutions
of elliptic equations was obtained independently by
Richard Casten and
Charles Holland, on the one side, and Hiroshi Matano,
on the other side, and it deals with the
classification of stable solutions with Neumann data. A paradigmatic
result in this case can be stated as follows:

\begin{theorem}[\cites{MR480282, MR555661}]\label{CH}
Let~$\Omega\subset\R^n$
be a bounded and convex domain with smooth boundary. 

Suppose that~$u$ is a smooth solution of
\begin{equation}\label{VAN00} \left\{
\begin{matrix}
-\Delta u(x)+ f(u(x))=0 & {\mbox{for any~$x\in\Omega$}}
\\
\displaystyle\frac{\partial u}{\partial\nu}(x)=0
& {\mbox{for any~$x\in\partial\Omega$,}}
\end{matrix}\right.
\end{equation}
for some smooth function~$f$, where~$\nu$ denotes the (external)
unit normal of~$\Omega$.

Assume also that~$u$ is stable, namely
\begin{equation}\label{VAN1}
\int_\Omega |\nabla\varphi(x)|^2+f'(u(x))\,|\varphi(x)|^2\,dx\ge0,
\end{equation}
for any~$\varphi\in H^1(\Omega)$.

Then, $u$ is necessarily constant.
\end{theorem}

We remark that 
\begin{equation}\label{VAN}
{\mbox{when $f$ vanishes identically
then~\eqref{VAN1} is automatically satisfied.}}\end{equation}
It is interesting to observe that,
with respect to Theorem~\ref{CH}, the fractional case
behaves very differently, and nonconstant stable solutions
with Neumann data
in convex domains do exist, according to the following result:

\begin{theorem}[\cite{SOAVE}]\label{CH2}
\label{ST} Let~$s\in(0,1)$.
There exist an open interval~$I\subset \R$
and a nonconstant function~$u$ such that
\begin{eqnarray}
\label{VAN2}
&& (-\Delta)^s u =0 \;{\mbox{ in }}\;I,\\
\label{VAN3}
&& \lim_{x\to x_0\in\partial I} \frac{u(x)-u(x_0)}{x-x_0}=0\\
\label{VAN4}
{\mbox{and }}&& u'=0 \;{\mbox{ on }}\;\partial I.\end{eqnarray} 
\end{theorem}

We observe that~\eqref{VAN2} is
the natural fractional counterpart
of~\eqref{VAN00} (with~$f:=0$, and~\eqref{VAN}
guarantees a stability condition). Of course, 
an interval is a (onedimensional) convex set, hence the
geometric setting of Theorem~\ref{CH}
is respected in Theorem~\ref{CH2}.
Also,
formula~\eqref{VAN4}
can be seen as a classical Neumann condition,
while
formula~\eqref{VAN3}
can be seen as a fractional Neumann condition
(say, of order~$s$).
Condition~\eqref{VAN3} is indeed quite exploited
as a natural boundary condition
in fractional problems, and it is compatible
with the boundary regularity theory and with the sliding methods,
see~\cite{MR3168912, MR3395749}
(for another notion of fractional Neumann condition
see~\cite{MR3651008}).

In this sense, Theorem~\ref{ST} can be considered as a ``counterexample''
for the fractional analogue of Theorem~\ref{CH} to hold.
The construction of Theorem~\ref{ST} is in fact very general. It is
based on Theorem~\ref{TH}
and provides a series of rather ``arbitrary'' counterexamples,
see Section~1.7 in~\cite{SOAVE} for additional details.

It has to be pointed out, however, that results similar
to the original ones in~\cites{MR480282, MR555661}
hold true for a different type of fractional
operator (the so-called ``spectral'' fractional Laplacian,
see~\cite{MR3233760}). In particular,
classification results for stable solutions of nonlocal operators
which can be seen as the fractional counterpart of those
in~\cites{MR480282, MR555661} have been given in
Sections~1.4, 1.5 and~1.6 in~\cite{SOAVE}.
This fact shows the very intriguing phenomenon,
according to which ``little'' modifications
in the fractional settings do produce rather different results,
which are sometimes in agreement with the classical case,
and sometimes not.

\subsection{The diffusive strategy of biological populations}

A classical problem in biomathematics
consists in studying the evolution of a biological species
with density~$u=u(x,t)$ in~$B_1\ni x$,
with prescribed boundary or external conditions.
In this framework, the so-called logistic equations
is based on the ansatz that the state of the population
is due to three well distinguishable features: 
\begin{itemize}
\item The population diffuses according to a stochastic motion;
\item For small density, the population grows more or less linearly,
thanks to some resources~$\rho=\rho(x)>0$;
\item When the density overcomes
a critical threshold~$\sigma/\mu$, for some~$\mu=\mu(x)>0$,
the population unfortunately
dies (roughly speaking, because ``there is no food for everybody'').
\end{itemize}
When the diffusion term is lead by the standard Brownian motion,
the logistic equation that we describe takes the form
\begin{equation}\label{LOGI}
\partial_t u=\Delta u + (\sigma-\mu u)\,u
\qquad{\mbox{ in }}B_1\times(0,T),\end{equation}
for some~$T>0$. In particular, the study of the steady states of~\eqref{LOGI}
leads to the equation
\begin{equation}\label{LOGI2}
-\Delta u = (\sigma-\mu u)\,u 
\qquad{\mbox{ in }}B_1.\end{equation}
On the other hand, recent experiments
have shown that several predators do not follow
standard diffusion processes, but rather discontinuous
processes with jumps whose distribution may exhibit a long (e.g. with a
polynomial tail), see e.g.~\cite{NATT}. This fact, that may seem
surprising, has indeed a sound motivation:
for a predator it makes little sense to move randomly looking
for prey, since, after a first attack, the other possible targets
will rapidly escape from the dangerous area --
conversely, a strategy of ``hit and run'',
based on quick hunts after long excursions, is more reasonable
to be efficient and ensure more food to the predator.

In this sense, a natural nonlocal variation of~\eqref{LOGI2} to be
taken into account is the fractional logistic equation
\begin{equation}\label{012938ersoer}
(-\Delta)^s u = (\sigma-\mu u)\,u 
\qquad{\mbox{ in }}B_1,\end{equation}
with~$s\in(0,1)$, see e.g.~\cites{MR3082317, MR3590678, MR3582231, MR3579567}
and the references therein.
Interestingly, different species in nature seem to exhibit different
values of the fractional parameter~$s$, probably due
to different environmental conditions and different morphological structures
and it is an intriguing problem to understand
what ``the optimal exponent~$s$'' should be in concrete
circumstances, see~\cite{MR3590646}.\medskip

Another interesting special feature offered by nonlocal
diffusion is the possibility for
nonlocal populations to efficiently plan their distribution
in order to consume
all (or almost all) the given resources in a certain ``strategic region''.
That is, if the region of interest for the population is, say, the ball~$B_1$,
the species can artificially and appropriately settle
its distribution outside~$B_1$, in order to satisfy in~$B_1$
a logistic equation as that in~\eqref{012938ersoer},
for a resource that is arbitrarily close to the original one.
The precise statement of this result is the following:

\begin{theorem}[\cites{MR3590678, MR3579567}]\label{BIO}
Assume that~$\sigma$, $\mu\in C^2(\overline{B_1})$, with
$$ \inf_{B_1}\sigma>0\quad{\mbox{ and }}\quad\inf_{B_1}\mu>0.$$
Then, for any~$\e>0$ there exist $u_\e$ and~$\sigma_\e$
such that
\begin{equation*}
\left\{ \begin{matrix}
\|\sigma-\sigma_\e\|_{C^2(B_1)}\le\e,\\
u_\e\ge \sigma_\e/\mu\quad{\mbox{ in }}B_1,\\
(-\Delta)^s u_\e =
(\sigma_\e-\mu u_\e)\,u_\e { \mbox{ in }}B_1
.\end{matrix}\right.
\end{equation*}
\end{theorem}

Once again, a proof of Theorem~\ref{BIO} may be
performed by exploiting Theorem~\ref{TH}, see Section~7
in~\cite{MR3579567}.

\section{Proof of Theorem~\ref{TH}}\label{ALL:FUN}

For simplicity, we focus
on the one-dimensional case: the general
case follows by technical modifications
and can be found in the original article~\cite{MR3626547}.

The core of the proof is to show that the derivatives of
$s$-harmonic functions have ``maximal span''
as a linear space
(and we stress that
this is not true for harmonic functions,
since the second derivatives of harmonic functions
satisfy a linear prescription).

We consider the set
\begin{equation}\label{calv} 
{\mathcal{V}}:= \big\{ h:\R\to\R
{\mbox{ s.t. $h$ is smooth and
$s$-harmonic in some neighborhood of the origin}} \big\}.\end{equation}
Notice that~${\mathcal{V}}$ has a linear space structure,
namely if~$h_1$
is $s$-harmonic in some open set~$V_1$ containing the origin
and~$h_2$
is $s$-harmonic in some open set~$V_2$ containing the origin, then,
for any~$\lambda_1$, $\lambda_2\in\R$, we have that~$h_3:=
\lambda_1 h_1+\lambda_2 h_2$
is $s$-harmonic in the open set~$V_3:=V_1\cap V_2\ni 0$.

Then, given~$J\in\N$, we define
\begin{equation}\label{calv2} {\mathcal{V}}_J := \big\{ \big(h(0),\,h’(0),\,\dots, \,h^{(J)}(0)\big)
{\mbox{ with }}h\in {\mathcal{V}}\big\}.\end{equation}
As customary, here~$h^{(J)}$ denotes the~$J$th derivative of the function~$h$.
In this way, we have that~${\mathcal{V}}_J$ is a linear subspace of~$\R^{J+1}$
(roughly speaking, each element of~${\mathcal{V}}_J$
is a~$(J+1)$-dimensional array containing the first~$J$ derivatives of a locally $s$-harmonic
function).

We claim that
\begin{equation}\label{TUTT}
{\mathcal{V}}_J=\R^{J+1}
\end{equation}
For this, we argue by contradiction
and we suppose that~${\mathcal{V}}_J$ is a linear subspace strictly
smaller than~$\R^{J+1}$. That is, ${\mathcal{V}}_J$ lies inside
a $J$-dimensional hyperplane, say with normal~$\nu$. Namely, there exists
\begin{equation}\label{TUTT2}\nu=(\nu_0,\dots,\nu_J)\in \R^{J+1}\;
{\mbox{ with }}\;|\nu|=1\end{equation}
such that
\begin{equation}\label{TUTT3}
{\mathcal{V}}_J\subseteq \big\{ X=(X_0,\dots,X_J)\in\R^{J+1} {\mbox{ s.t. }}\nu\cdot X=0
\big\} 
\end{equation}
Now, for any~$ { t }>0$, we define
$$h_{ t } (x):= (x+ { t })_+^s.$$
It is known that~$h_{ t } $ is~$s$-harmonic
in~$(- { t },+\infty)$ (see e.g. Chapter~2.4 in~\cite{MR3469920}
for an elementary proof). Consequently,
$h_{ t } \in{\mathcal{V}}$ and then
$$ X _{ t } :=\big(h_{ t } (0),\,\dots,\,h^{(J)}_{ t } (0)\big)
\in{\mathcal{V}}_J.$$
As a result, by~\eqref{TUTT3},
\begin{equation}\label{TUTTX} 0= \nu\cdot X _{ t } =\sum_{j=0}^J
\nu_j h_{ t }^{(j)}(0)=\sum_{j=0}^J
\mu_{s,j}\, { t }^{s-j},\end{equation}
where
\begin{equation}\label{CAN} \mu_{s,j}:= \nu_j \prod_{i=0}^{j-1} (s-i).\end{equation}
Hence, multiplying the identity in~\eqref{TUTTX}
by~$t^{J-s}$,
for any~$t>0$, it holds that
$$ \sum_{k=0}^J
\mu_{s,J-k}\, { t }^{k}=0,$$
which, by the
Identity Principle for Polynomials,
implies that~$\mu_{s,0}=\dots=\mu_{s,J}=0$ and accordingly\footnote{We
stress that here it is crucially used the fact that~$s$
is not an integer.}
from~\eqref{CAN} we get that~$\nu_{0}=
\dots=\nu_{J}=0$.
This is in contradiction with~\eqref{TUTT2}
and so the proof of~\eqref{TUTT}
is complete.

Now, the proof of Theorem~\ref{TH} follows by approximation and scaling.
Given~$\bar v\in C^2(\overline{B_1})$
and~$\e\in(0,1)$, 
in view of the Stone-Weierstrass Theorem
we take a polynomial~${{\mathcal{P}}_\e}$ such that
\begin{equation}\label{PESP}
\|\bar v-{{\mathcal{P}}_\e}\|_{C^2(B_1)}\le\frac\e2.
\end{equation}
We write
$$ {{\mathcal{P}}_\e}(x)=\sum_{j=0}^{N_\e} c_{j,\e} x^j=\sum_{j=0}^{N_\e} m_{j,\e}(x),$$
for some~$N_\e\in\N$ and some~$c_{1,\e},\dots,c_{N_\e,\e}\in \R$,
where
\begin{equation}\label{A24}
m_{j,\e}(x):=c_{j,\e} x^j.
\end{equation}
Without loss of generality, by possibly adding zero coefficients in the representation above,
we can suppose that
\begin{equation}\label{N:3}
N_\e\ge3.\end{equation}
We set
$$C_\e:=\max_{j\in\{0,\dots,N_\e\}} |c_{j,\e}|.$$
For any~$j\in\{0,\dots,N_\e\}$, we let~$H_{j,\e}:\R\to\R$
be a function which is $s$-harmonic in a neighborhood of the origin and such that,
for any~$i\in\{0,\dots,N_\e\}$ it holds that
\begin{equation}\label{A25} H_{j,\e}^{(i)}(0)=\left\{
\begin{matrix}
c_{j,\e} \;{ j! } &{\mbox{ if }} i=j,\\
0 &{\mbox{ otherwise.}}
\end{matrix}
\right.\end{equation}
Once again, $H_{j,\e}^{(i)}$ denotes here the~$i$th derivative
of~$H_{j,\e}$.
We stress that the existence of~$H_{j,\e}$ is a consequence
of~\eqref{TUTT}.
We also set
\begin{equation}\label{A26}\begin{split}
& r_{j,\e}:= \frac{\e}{10\,N_\e^2\,\Big(1+
\displaystyle\sup_{x\in(-1,1)}\big|H_{j,\e}^{(N_\e+1)}(x)\big| \Big)}\in(0,1)
\\{\mbox{ and }}\;&
{\mathcal{H}}_{j,\e}(x):= r_{j,\e}^{-j}\; H_{j,\e}(r_{j,\e} x). \end{split}\end{equation}
We remark that,
for any~$i$, $j\in\{0,\dots,N_\e\}$,
\[ {\mathcal{H}}_{j,\e}^{(i)}(0)=\left\{
\begin{matrix}
c_{j,\e} \;{ j! } &{\mbox{ if }} i=j,\\
0 &{\mbox{ otherwise,}}
\end{matrix}
\right.\]
thanks to~\eqref{A25}.
Therefore, in view of~\eqref{A24}, the function
\begin{equation}
\label{A28}
{\mathcal{D}}_{j,\e}(x):=
{\mathcal{H}}_{j,\e}(x)-c_{j,\e} x^j=
{\mathcal{H}}_{j,\e}(x)-m_{j,\e}(x)\end{equation}
satisfies
\begin{equation}\label{Dj}
{\mathcal{D}}_{j,\e}^{(i)}(0)=0\;{\mbox{ for all }}\; i\in\{0,\dots,N_\e\}.
\end{equation}
In addition, for any~$x\in(-1,1)$ and any~$j\in\{0,\dots,N_\e\}$,
\begin{eqnarray*}
\big| {\mathcal{D}}_{j,\e}^{(N_\e+1)}(x)\big| &=&
\big|{\mathcal{H}}_{j,\e}^{(N_\e+1)}(x)\big|\\
&\le&
r_{j,\e}^{N_\e+1-j} \big|H_{j,\e}^{(N_\e+1)}(r_{j,\e} x)\big|
\\ &\le& 
r_{j,\e}\sup_{(-1,1)}\big|H_{j,\e}^{(N_\e+1)}\big|
\\ &\le& \frac{\e}{2 N_\e^2},
\end{eqnarray*}
thanks to~\eqref{A26}.
This, \eqref{Dj} and a Taylor expansion 
give that, for any~$x\in(-1,1)$
and any~$i$, $j\in\{0,\dots,N_\e\}$,
$$ \big| {\mathcal{D}}_{j,\e}^{(i)}(x)\big|\le
\sup_{(-1,1)}
\big| {\mathcal{D}}_{j,\e}^{(N_\e+1)}\big|
\le \frac{\e}{10\,N_\e^2}.$$
Hence, recalling~\eqref{N:3}
$$ \sum_{j=0}^{N_\e} \| {\mathcal{D}}_{j,\e}\|_{C^2(-1,1)}
\le \frac{\e}{2}.$$
So, we define
$$ v_\e:=\sum_{j=0}^{N_\e}{\mathcal{H}}_{j,\e}.$$
We have that~$v_\e$ is $s$-harmonic in~$(-1,1)$ and, recalling~\eqref{PESP}
and~\eqref{A28},
\begin{eqnarray*}
\|\bar v-v_\e\|_{C^2(-1,1)} &\le&
\|\bar v-{{\mathcal{P}}_\e}\|_{C^2(-1,1)}+
\|{{\mathcal{P}}_\e}-v_\e\|_{C^2(-1,1)}\\
&\le& \frac\e2 + \left\|\sum_{j=0}^{N_\e}\left(m_{j,\e}-  
{\mathcal{H}}_{j,\e}\right)\right\|_{C^2(-1,1)}
\\&\le& \frac\e2 + \sum_{j=0}^{N_\e}\left\| 
{\mathcal{D}}_{j,\e}\right\|_{C^2(-1,1)}
\\ &\le&\e.
\end{eqnarray*}
This establishes Theorem~\ref{TH} in this setting.

\vfill

\end{document}